\documentclass{amsart}
 \usepackage{amssymb,amsmath,amsfonts,epsfig,latexsym}

 \newtheorem{theorem}{Theorem}[section]
\newtheorem{definition}[theorem]{Definition}
\newtheorem{proposition}[theorem]{Proposition}
\newtheorem{lemma}[theorem]{Lemma}

\newtheorem*{EKLT}{Einsiedler-Kapranov-Lind Theorem}
\newtheorem*{ZT}{Zhang's Theorem}

\theoremstyle{definition}
\newtheorem{remark}[theorem]{Remark}
\newtheorem{example}[theorem]{Example}

\def\Z{\ensuremath{\mathbb{Z}}}

\def\C{\ensuremath{\mathbb{C}}}
\def\R{\ensuremath{\mathbb{R}}}
\def\G{\ensuremath{\mathbb{G}}}

\def\A{\ensuremath{\mathcal{A}}}

\def\adele{\A_\mathbb{A}}

\def\<{\ensuremath{\langle}}
\def\>{\ensuremath{\rangle}}

\DeclareMathOperator{\ev}{ev}

\DeclareMathOperator{\Hom}{Hom}

\DeclareMathOperator{\relint}{int}

\begin{document}

\title[Adelic amoebas]{Adelic amoebas disjoint from open halfspaces}

\author[Payne]{Sam Payne}

\thanks{This work was done during a visit to the Centre de Recherches Math\'ematiques in Montr\'eal, Canada.  Supported by the Clay Mathematics Institute.}

\begin{abstract}
We show that a conjecture of Einsiedler, Kapranov, and Lind on adelic amoebas of subvarieties of tori and their intersections with open halfspaces of complementary dimension is false for subvarieties of codimension greater than one that have degenerate projections to smaller dimensional tori.  We prove a suitably modified version of the conjecture using algebraic methods, functoriality of tropicalization, and a theorem of Zhang on torsion points in subvarieties of tori.
\end{abstract}

\maketitle

\section{Introduction}

Let $X$ be a subvariety of a torus over a field $K$ that is either a number field or the function field of a curve.  The adelic amoeba of $X$ is a union of amoebas $\A_p(X)$, one for each place $p$ of $K$, and the intersections of adelic amoebas with rational open halfspaces govern nonexpansive sets in algebraic dynamical systems.  See \cite{Schmidt90}, \cite{Schmidt95}, and \cite[Section~7]{BoyleLind97} for background and further references on algebraic dynamical systems and expansive subdynamics in general, and \cite{ELMW} and \cite[Section~4]{EKL} for details on the relationship between expansive subdynamics and amoebas. 

\begin{EKLT} \cite[Theorem~2.3.3]{EKL} Let $X$ be a hypersurface in a torus over a number field or the function field of a curve.  If there is an open half line that is disjoint from the adelic amoeba of $X$ then $\A_p(X)$ contains zero for every $p$.
\end{EKLT}

\noindent A special case of the Einsiedler-Kapranov-Lind Theorem was proved using homoclinic points for algebraic dynamical systems in \cite[Proposition~5.5]{ELMW}.  The full result was proved using analytic methods involving domains of convergence of Laurent series expansions of the reciprocal of the equation defining $X$ with respect to both archimedean and nonarchimedean norms.  

Einsiedler, Kapranov, and Lind then conjectured that a similar result should hold for adelic amoebas of higher codimension subvarieties. Specifically, they conjectured that if $X$ is $r$-dimensional and the adelic amoeba $\adele(X)$ is disjoint from some codimension $r$ halfspace, then either $\A_p(X)$ contains zero for every $p$ or $\adele(X)$ is contained in a hyperplane \cite[Conjecture~2.3.5]{EKL}.  However, this conjecture is false in general, for both function fields and number fields, when there is a projection, or split surjection of tori, whose restriction to $X$ is degenerate in the sense that the image of $X$ either has smaller than expected dimension or is cut out by an equation whose coefficients lie in the subfield $k$ of a function field $k(C)$.  See Examples~\ref{curve in 3space} and \ref{surface in 4space}, below.  Roughly speaking, our main result says that if the adelic amoeba is not contained in a hyperplane then these are the only obstructions to the adelic amoeba meeting any given rational open halfspace of complementary dimension.  

Let $T$ be the torus with character lattice $M$, and let $N = \Hom(M,\Z)$ be the dual lattice, with $N_\R$ the real vector space $N \otimes \R$.  An open halfspace $H$ in $N_\R$ is the sum of a linear subspace $\partial H$ and an open halfline $\R_{>0} \cdot v$ disjoint from $\partial H$.  When $\partial H$ is rational we write $T_{\partial H}$ for the subtorus of $T$ whose lattice of one parameter subgroups is $\partial H \cap N$ and $\varphi$ for the projection from $T$ to $T/ T_{\partial H}$.  Let $X'$ be the closure of the image of $X$ in $T/ T_{\partial H}$.

\vspace{5 pt}

\begin{theorem} \label{main}
Let $X$ be a $r$-dimensional subvariety of $T$ defined over a global field $K$.  If $H$ is disjoint from the adelic amoeba of $X$ then either
\begin{enumerate}
\item $X'$ has codimension greater than one, 
\item $K$ is a function field $k(C)$ and $X'$ is defined over $k$, or
\item $K$ is a number field and $X'$ is the translate of a subtorus by a torsion point.
\end{enumerate}
\end{theorem}

\noindent In the special case where $X$ is a hypersurface, $\varphi$ is just the identity on $T$.  If $K$ is a function field and $X$ is defined over $k$ then $\A_p(X)$ is the normal fan of the Newton polytope of a defining equation, and hence contains zero, for every $p$. Similarly, if $K$ is a number field and $X$ is a translate of a subtorus by a torsion point, then $\A_p(X)$ is a hyperplane, and hence contains zero, for every $p$.  In particular, we recover the Einsiedler-Kapranov-Lind Theorem.  Furthermore, Theorem~\ref{main} is sharp in the following sense.  If $K$ is a function field and there is a projection from $T$ to an $r+1$-dimensional quotient torus such that the closure $X'$ of the image of $X$ either has codimension greater than one or is defined over $k$, then a general open half line is disjoint from the adelic amoeba of $X'$, and its preimage is an open halfspace disjoint from the adelic amoeba of $X$.  If $K$ is a number field and $X'$ is a translate of a subtorus by a torsion point then the adelic amoeba of $X$ is contained in a hyperplane.

To prove Theorem~\ref{main}, we first treat the hypersurface case in Propositions~\ref{hypersurface:ff} and \ref{hypersurface:nf}.  For function fields, the proof is self-contained and algebraic.  For number fields, we obtain best possible results using Zhang's Theorem on the Zariski closure of the set of torsion points in a subvariety of an algebraic torus \cite{Zhang95}.  See Section~\ref{section:nf} for details.  The full result then follows from the hypersurface case by functoriality of tropicalization for the projection $\varphi$, and an analogue of the density of fibers of tropicalization for archimedean amoebas (Lemma~\ref{archimedean density}). 

\vspace{5 pt}

We conclude the introduction with two counterexamples to the conjecture mentioned above, where $H$ is a rational halfspace disjoint from $\adele(X)$, but $\A_p(X)$ does not contain zero for some place $p$, and $\adele(X)$ is not contained in a hyperplane.  In the first example, $X'$  is a hypersurface defined over the scalar subfield of a function field.

\begin{example} \label{curve in 3space}
Suppose $K$ is the function field $\C(z)$.  Let $X$ be the image of $\G_m \smallsetminus \{ 1,  z^{-1} \}$ in $\G_m^3$ under the map
\[
t \mapsto (t, t-1, t - z^{-1}).
\]
It is straightforward to check that at all places $p$ of $K$ other than $z = 0$ and $z = \infty$, the amoeba $\A_p(X)$ is the union of four rays in $\R^3$, spanned by $e_1$, $e_2$, $e_3$, and $-e_1 -e_2 -e_3$, respectively.  In particular, $\adele(X)$ is not contained in any hyperplane.  Furthermore, if $t$ is regular and nonvanishing at $z = 0$, then $t - z^{-1}$ has a pole at zero.  It follows that zero is not contained in the amoeba of $X$ at the place $z = 0$.  However, the adelic amoeba $\adele(X)$ does not meet the open halfspace $\R_{>0} (e_1 + e_2) + \R \cdot e_3$, since $t$ and $t-1$ cannot both vanish at the same point.
\end{example}

\noindent In the next example, the projected image $X'$ is not a hypersurface.

\begin{example} \label{surface in 4space}
Suppose $K$ is the field of rational numbers.  Let $X$ be the surface in $\G_m^4$ that is the image of $(\G_m \smallsetminus \{ 1, 2 \}) \times \G_m$ under the map
\[
(t,t') \mapsto (t, t-1, t-2, t').
\]
The amoeba of $X$ at the archimedean place, which is the usual complex amoeba studied in \cite{GKZ}, does not contain zero, since complex numbers $t$, $t-1$, and $t -2$ cannot all lie on the unit circle.  Furthermore, it is straightforward to check that at every nonarchimedean place except two, $\A_p(X)$ contains the rays spanned by $e_1$, $e_2$, $e_3$, and $e_4$, so $\adele(X)$ is not contained in a hyperplane.  Let $H$ be the halfspace
\[
H = \R_{> 0} \cdot (e_1 + e_2 + e_3) + \R \cdot e_4
\]
in $\R^4$.  We claim that $H$ does not meet $\adele(X)$.  This is clear because, at an archimedean place, complex numbers $t$, $t-1$, and $t-2$ cannot all have the same absolute value, and at a nonarchimedean place, $t$ and $t-1$ cannot both have positive valuation.
\end{example}

\noindent It is straightforward to construct examples similar to Example~\ref{surface in 4space} over function fields.

\vspace{5 pt}

\noindent \textbf{Acknowledgments.}  I am grateful to M. Einsiedler, D. Savitt, and R. Vakil for helpful discussions.

\section{Preliminaries}

Let $K$ be either a number field or the function field of a curve, and let $S$ be the set of places of $K$.  Let $| \ |_p$ denote the normalized absolute value representing a place $p \in S$.  Let $K_p$ be the completion of $K$ with respect to $| \ |_p$.  Recall that $| \ |_p$ extends uniquely to the algebraic closure $\overline K_p$  by setting $|x|_p = | N_{L/K_p}(x)|^{1/[L:K_p]}$, for $x$ contained in a finite extension $L/K_p$.  Furthermore, we have the product formula
\[
\prod_{p \in S} | a |_p = 1,
\]
for all $a \in K ^*$ \cite[Chapter~2]{CasselsFrohlich67}.

Let $T$ be a torus with character lattice $M$, and let $N_\R = \Hom (M, \R)$.  A point $x \in T(\overline K_p)$ corresponds to an evaluation map
\[
\ev_x : M \rightarrow \overline K_p^*,
\]
given by $u \mapsto \chi^u(x)$.  For each $p$, composing with $-\log | \ |_p$ gives a group homomorphism from $M$ to $\R$, and hence a point $-\log| x |_p \in N_\R$.  If we choose compatible coordinates $T \cong \G_m^n$ and $N_\R \cong \R^n$, then $x$ is given by a tuple $(x_1, \ldots, x_n)$ of nonzero elements of $\overline K_p$, and $-\log| x |_p$ is the vector
\[
-\log|x|_p = (-\log| x_1|_p , \ldots, -\log| x_n |_p).
\]
For a subvariety $X$ in $T$ defined over $K$, the amoeba of $X$ at $p$ is the closure of the image of $X(\overline K_p)$ in $N_\R$,
\[
\A_p(X) = \overline{ \{ -\log| x |_p \ : \ x \in X(\overline K_p) \} }.
\]

\begin{remark} 
If $K$ is a number field and $p$ is a finite place, or if $K$ is a function field, then $-\log| \ |_p$ is a nonarchimedean valuation and $\A_p(X)$ is the underlying set of a polyhedral complex of pure dimension equal to the dimension of $X$, and is the closure of the ``tropicalization" of $X$ with respect to this nonarchimedean valuation.  Furthermore, if $G \subset \R$ is the valuation group, then the image of $X(\overline K)$ is the set of $G$-rational points in $\A_p(X)$, and the preimage of any $G$-rational point in $\A_p(X)$ is Zariski dense in $X(\overline K)$, by  \cite[Theorem~4.1]{tropicalfibers}.  If $K$ is a number field and $\infty$ is an infinite place corresponding to an embedding of $K$ in $\C$, then $\A_\infty(X)$ is the usual amoeba of the complex variety $X \times_K \C$, studied in \cite{GKZ}.
\end{remark}

\begin{definition} The adelic amoeba $\adele(X)$ is the union
\[
\adele(X) = \bigcup_{p \in S} \A_p(X).
\]
\end{definition}

\section{Adelic amoebas of hypersurfaces over function fields}

Let $K = k(C)$ be the function field of a curve.  Let $f = a_1 x^{u_1} + \cdots + a_s x^{u_s}$ be a Laurent polynomial in $K[M]$, with $a_i \in K^*$ and $s \geq 2$.  Let $Q$ be the Newton polytope of $f$, and let $X = V(f)$ be the hypersurface in $T$ cut out by $f$.

For any place $p$ of $K$, the nonarchimedean amoeba $\A_p(X)$ is the corner locus of the convex piecewise linear function on $N_\R$ given by
\[
\Psi_f(v) = \min \{ \<u_i, v\> + \nu_p(a_i) \},
\]
where $\nu_p$ is the valuation at $p$.  In particular, for all but finitely many $p$, $\nu_p(a_i)$ is zero and $\A_p(X)$ is the union of the codimension one cones in the inward normal fan $\Delta$ of $Q$.  Therefore, if there is an open half line $\R_{>0} \cdot v$ in $N_\R$ that is disjoint from $\adele(X)$, then it must lie in the interior of some maximal cone of $\Delta$.  The maximal cones $\sigma_i$ of $\Delta$ correspond to those $u_i$ that are vertices of $Q$, and the interior of $\sigma_i$ is
\[
\relint(\sigma_i) = \{ v \in N_\R : \<u_i, v\> > \< u_j, v\> \mbox{ for all } j \neq i \}.
\]

\begin{lemma} \label{unbounded components}
An open half line in $\relint(\sigma_i)$ is disjoint from $\A_p(X)$ if and only if $\nu_p(a_i)$ is less than or equal to $\nu_p(a_j)$ for all $j$.
\end{lemma}

\begin{proof}
Let $\R_{>0} \cdot v$ be an open half line in $\relint(\sigma_i)$.  Let $\Psi_j$ be the affine linear function on $\R$ given by
\[
\Psi_j(c) = \<u_j, cv\> + \nu_p(a_j) - \big( \<u_i, cv\> + \nu_p(a_i) \big) .
\]
Since $v$ is in the interior of $\sigma_i$, $\Psi_j(c)$ is positive for $c \gg 0$.  If $\nu_p(a_j)$ is greater than $\nu_p(a_i)$ for some $j$, then $\Psi_j(0)$ is less than zero and there is a unique positive real number $c_j$ such that $\Psi_j(c_j) = 0$.  If $c_k$ is the largest such number, then $c_k \cdot v$ is in $\A_p(X)$.

Conversely, if $\nu_p(a_i)$ is less than or equal to $\nu_p(a_j)$ for all $j$, then $\Psi_j(c \cdot v)$ is positive for all $c > 0$, and it follows that the open half line $\R_{>0} \cdot v$ is disjoint from $\A_p(X)$.
\end{proof}

\begin{proposition} \label{hypersurface:ff}
Let $X$ be a hypersurface in $T$ defined over a function field $k(C)$.  Then $\adele(X)$ is disjoint from some open half line in $N_\R$ if and only if $X$ is defined over $k$.
\end{proposition}

\begin{proof}
If $X$ is defined over $k$, then every amoeba $\A_p(X)$ is equal to the codimension one skeleton of the normal fan of the Newton polytope of a defining equation for $X$, and hence any open halfline that is contained in a maximal cone of this fan is disjoint from $\adele(X)$.

Suppose $\R_{>0} \cdot v$ is an open half line in $N_\R$ that is disjoint from $\adele(X)$.  Then $v$ lies in $\relint(\sigma_i)$ for some $i$.  By Lemma~\ref{unbounded components}, $\nu_p(a_i)$ is less than or equal to $\nu_p(a_j)$ for all $p$ and all $j$.  Then, for every $j$, the quotient $a_j/a_i$ is a rational function on $C$ with no poles, and hence lies in $k$.  Therefore, $f/a_i$ is a defining equation for $X$ with coefficients in $k$.
\end{proof}

\noindent For hypersurfaces in tori over function fields, Proposition~\ref{hypersurface:ff} is stronger than the Einsiedler-Kapranov-Lind Theorem, since there do exist hypersurfaces that are not defined over $k$ such that $\A_p(X)$ contains zero for every $p$.  The following example illustrates this possiblity.

\begin{example}
Suppose $K$ is the function field $\C(z)$.  Let $X$ be the curve in $\G_m^2$ defined by the equation
\[
z x_1 + (z-1) x_2 + (z-2) = 0.
\]
For all places $p$ other than $z =0$, $z = 1$, and $z = 2$, $\A = \A_p(X)$ is the union of the rays spanned by $e_1$, $e_2$, and $-e_1-e_2$.  At the places $z = 0$, $z = 1$, and $z = 2$, the amoebas of $X$ are $\A - e_1$, $\A-e_2$, and $\A + e_1 + e_2$, respectively, all of which contain zero.
\end{example}

\section{Adelic amoebas of hypersurfaces over number fields} \label{section:nf}

Over number fields, one of our main tools is the characterization of torsion points in terms of normalized absolute values.

\begin{lemma}  \label{torsion lemma}
Let $K$ be a number field, and let $x$ be a point in $T(\overline K)$.  If $-\log |x|_p$ and $-\log |x|_\infty$ are equal to zero for every place $p$ of $K$ and every infinite place $\infty$ of a finite extension $L/K$ over which $x$ is defined, then $x$ is a torsion point.
\end{lemma}

\begin{proof}
Choose coordinates $T \cong (\G_m)^n$, so $x$ is given by a tuple $(x_1, \ldots, x_n)$ of nonzero elements of $L$.  If $-\log|x|_p = 0$ then $|x_i|_p= 1$ for each place $p$ of $K$.  Then $x_i$ is a unit in the ring of integers of $L$, so $|x_i|_q = 1$ for every finite place $q$ of $L$.  Since $|x_i|_\infty = 1$ for every infinite place of $L$, by hypothesis, it follows that $x_i$ is a root of unity \cite[p.~72]{CasselsFrohlich67}.
\end{proof}

\noindent This characterization of torsion points is especially useful because of the following analogue of Bogomolov's Conjecture for algebraic tori, due to Zhang  \cite{Zhang95}.  See also the elementary proof of Bombieri and Zannier \cite{BombieriZannier95}.

\begin{ZT}
Let $X$ be a subvariety of $T$ defined over a number field.  Then the Zariski closure of the set of torsion points in $X(\overline K)$ is a finite union of translates of subtori by torsion points.
\end{ZT}

\noindent  Zhang's Theorem leads to the following characterization of hypersurfaces in $T$ whose adelic amoebas are disjoint from some open half line.

\begin{proposition} \label{hypersurface:nf}
Let $X$ be a hypersurface in $T$ defined over a number field.  Then $\adele(X)$ is disjoint from some open half line if and only if $X$ is a translate of a subtorus by a torsion point.
\end{proposition}

\begin{proof}
If $X$ is a translate of a subtorus by a torsion point, then $\adele(X)$ is a hyperplane, and hence is disjoint from any half line that does not lie in that hyperplane.

Suppose $\adele(X)$ is disjoint from some half line $\R_{>0} \cdot v$.  Choose a complete fan $\Delta$ that refines the normal fan of the Newton polytope of a defining equation for $X$.  Since the set of open half lines disjoint from $\adele(X)$ is open, we may assume that $v$ is rational and the line $\R \cdot v$ meets the codimension one skeleton of $\Delta$ only at zero.  Let $T_0$ be the one parameter subgroup corresonding to $\R \cdot v$, and let $t \in T(\overline K)$ be a torsion point.  

We claim that $t T_0 \cap X$ is nonempty.  To see this, note that the closure of $X$ in the toric  variety associated to $\Delta$ is Cartier and globally generated and defines a morphism that does not contract the closure of $t T_0$, so the intersection number $(\overline X \cdot \overline {tT_0})$ is positive.  Furthermore, $\overline X$ does not contain any $T$-fixed points, since $\A_p(X)$ does not intersect the relative interiors of any of the maximal cones of $\Delta$, and in particular $\overline X$ does not contain either of the points in $\overline {tT_0} \smallsetminus tT_0$.  This property of nonarchimedean amoebas is proved in \cite[Lemma~2.2]{Tevelev07} over fields of Puiseux series; the proof for number fields with respect to a $p$-adic valuation is similar.  Hence $tT_0$ intersects $X$ in finitely many points.

Next, we claim that any point $x$ in $t T_0 \cap X$ is torsion.  By Lemma~\ref{torsion lemma}, it suffices to show that $-\log|x|_p$ and $-\log|x|_\infty$ are equal to zero for all places $p$ of $K$ and all infinite places of some extension $L/K$ over which $x$ and $t$ are both defined.  For places $p$ of $K$, this is clear, since $-\log |x|_p$ lies on $\A_p(tT_0) = \R_v$ for all $v$, and the product formula implies that $\sum_p -\log|x|_p = 0$.  Since the open half line $\R_{>0} \cdot v$ is disjoint from $\adele(X)$, by hypothesis, it follows that $-\log|x|_p = 0$ for all $p$.  For the infinite places, if $-\log|x|_\infty$ is nonzero, then the product formula for $L$ implies that there is a Galois conjugate complex point $x'$, which lies in the intersection of $X$ with a translation of $T_0$ by a conjugate torsion point $t'$, such that $-\log|x'|_\infty$ lies on $\R_{>0} \cdot v$, contradicting the hypothesis that $\adele(X)$ is disjoint from $\R_{>0} \cdot v$.

Now we have shown that for any torsion point $t$, $X \cap tT_0$ contains a torsion point.  If the torsion points in $X$ were contained in a finite union of translates of subtori of codimension greater than one in $T$, then we could choose a torsion point $t$ such that $tT_0$ contained no torsion points of $X$.  Therefore, by Zhang's Theorem it follows that $X$ is a translate of a codimension one subtorus by a torsion point.
\end{proof}

\noindent For hypersurfaces in tori over number fields, Proposition~\ref{hypersurface:nf} is stronger than the Einsiedler-Kapranov-Lind Theorem, since there do exist hypersurfaces that are not translates of subtori by torsion points such that $\A_p(X)$ contains zero for every $p$.  The following example illustrates this possiblity.

\begin{example}
Suppose $K$ is the field of rational numbers.  Let $X$ be the curve in $\G_m^2$ defined by the equation
\[
x_1x_2 -2x_1 -2x_2 + 1 = 0.
\]
For all places $p$ other than $p = 2$ and $p = \infty$, $\A_p(X)$ is the union of the coordinate axes.  Let $v = (1,-1)$.  Then $\A_2(X)$ is the union of the segment $(v,-v)$ and the rays based at $\pm v$ in directions $\pm e_1$ and $\mp e_2$, respectively.  At the infinite place, $\A_\infty(X)$ has a ``pinching point," in the sense of \cite[Section~3.5.3]{Mikhalkin04}, at zero.
\end{example}

\section{Proof of main result}

Let $\phi: N_\R \rightarrow N'_\R$ be the map of vector spaces induced by $\varphi$, where $N'$ is the character lattice of $T/ T_{\partial H}$.  We will deduce Theorem~\ref{main} from Propositions~\ref{hypersurface:ff} and \ref{hypersurface:nf}, by comparing the adelic amoeba of $X'$ with the image of $\adele(X)$ under $\phi$.  For any place $p$, $\phi$ restricts to give a map
\[
\phi:  \A_p(X) \rightarrow \A_p(X'),
\]
by functoriality of tropicalization \cite[Section~2]{tropicalfibers}.  If $p$ is a nonarchimedean place, then this map is surjective.  This surjectivity was proved by Tevelev \cite[Proposition~3.1]{Tevelev07} for fields of Puiseux series and may be seen as an immediate consequence of the Zariski density of the fibers of tropicalization \cite[Corollary~4.2]{tropicalfibers}.  At an archimedean place, the natural map $\phi: \A_\infty(X) \rightarrow \A_\infty(X')$ is not surjective in general, as the following example shows.

\begin{example}
Let $X$ be the image of $\G_m \smallsetminus \{ 1, 2 \}$ in $\G_m^3$ under the map
\[
t \mapsto (t, \, 2-t, \, t-1).
\]
The projection to the first two factors maps $X$ into the curve $X'$ in $\G_m^2$ where $x_1 + x_2 = 2$, and the image of $X$ is $X' \smallsetminus (1,1)$.  Since $(1,1)$ is the unique point in $X'$ with both coordinates on the unit circle, it follows that the image of $\A_\infty(X)$ is exactly $\A_\infty(X') \smallsetminus 0$.
\end{example}

\noindent To overcome this minor complication at the archimedean places, we use the following lemma.

\begin{lemma} \label{archimedean density}
Let $X$ be a positive dimensional hypersurface in $T$.  Then the preimage of any rational open interval in $N_\R$ that meets $\A_\infty(X)$  is Zariski dense in $X(\overline K_\infty)$.
\end{lemma}

\begin{proof}
Fix a two dimensional subtorus $T_0 \subset T$.  Then $X$ is the union of its intersections with the translates of $T_0$, and the intersection of $T$ with a general translate of $T_0$ is an irreducible curve.  Therefore, it will suffice to prove the lemma in the case where $T$ is two-dimensional.

Suppose $T$ is two-dimensional, $X$ is an irreducible curve in $T$, and $(v_1, v_2)$ is a rational interval that intersects $\A_\infty(X)$. Choose a splitting $T \cong T_1 \times T_2$, where $T_1$ is the one parameter subgroup corresponding to $\R \cdot (v_1 - v_2)$.  Then either $X$ is a translate of $T_1$, in which case the lemma is obvious, or there is an isolated point $x$ in the intersection of $X$ with the translate of $T_1$ by a point $t$ in $T_2$ and $-\log| x |_\infty$ lies in $(v_1, v_2)$.  Then, since $-\log| \ |_\infty$ is continuous, for torsion points $t'$ sufficiently closet to the identity in $T_2$, the translation of $T_1$ by $t't$ intersects $X$ in a point sufficiently close to $x$ so that its image under $-\log| \ |_\infty$ lies in the interval $(v_1, v_2)$.  In particular, the set of points in $X(\overline K_\infty)$ in the preimage of $(v_1, v_2)$ is infinite, and hence Zariski dense.
\end{proof}

\begin{proof}[Proof of Theorem~\ref{main}]
Assume $X'$ is a hypersurface and $H$ is disjoint from $\adele(X)$.  If $K$ is a function field then the adelic amoeba is the image of $\adele(X)$ under $\phi$.  Then the open half line $\phi(H)$ is disjoint from $\adele(X')$, and $X'$ is defined over $k$, by Proposition~\ref{hypersurface:ff}.

Let $K$ be a number field.  We claim that the open half line $\phi(H)$ is disjoint from $\adele(X')$.  At a finite place $p$, $\A_p(X')$ is the image of $\A_p(X)$, so $\phi(H)$ is disjoint from $\A_p(X')$. At an infinite place, by Lemma~\ref{archimedean density} and the fact that $\varphi(X)$ contains an open dense subvariety of $X'$, if $\phi(H)$ meets $\A_\infty(X')$ then $H$ meets $\A_\infty(X)$.  We conclude that the open half line $\phi(H)$ is disjoint from $\adele(X')$, and hence $X'$ is a translate of a subtorus by a torsion point, by Proposition~\ref{hypersurface:nf}.
\end{proof}

\bibliography{math}
\bibliographystyle{amsalpha}

\end{document}